\input amstex
\documentstyle{amsppt}
\input bull-ppt
\keyedby{bull261e/PAZ}
\define\Bysame{\hbox to 2em{\hrulefill}\thinspace,}
\define\REF#1.#2\par {[#1]\quad#2\par} 
\define\nr{\Vert}
\define\litspace{\mskip\medmuskip}
\define\bfa{\bold{a}}
\define\bfb{\bold{b}}
\define\bfc{\bold{c}}
\define\cvt{2^{\aleph_0}}
\define\cardset{(\lambda_i)_{i\in I}}
\define\includedin{\subseteq}
\define\intersect{\cap}
\define\Aa{{\aleph_\alpha}}
\define\Ad{{\aleph_\delta}}
\define\Ao{{\aleph_\omega}}
\define\Az{{\aleph_0}}

\define\scriptbb#1{\rlap{$\rm \scriptstyle 
#1$}\kern.5pt\hbox{$\rm
 \scriptstyle#1$}}

\define\littleslant{\hbox{$\scriptscriptstyle/$}}
\define\slantline{\raise 1.5 pt
\hbox{\rlap{\littleslant}\kern.1pt\raise.2pt\littleslant}}
\define\add#1#2#3{\rlap{\kern #3 pt #1}{\rm #2}}

\define\FF{{\cal F}}
\define\JJ{{\cal J}}

\define\cl{\operatorname{cl}}
\define\cof{\operatorname{cof}}
\define\cov{\operatorname{cov}}
\define\pcf{\operatorname{pcf}}
\define\pp{\operatorname{pp}}
\define\tcf{\operatorname{tcf}}
\define\PP{\operatorname{PP}}
\define\cf{\operatorname{cf}}
\topmatter
\cvol{26}
\cvolyear{1992}
\cmonth{April}
\cyear{1992}
\cvolno{2}
\cpgs{197-210}

\title Cardinal arithmetic for skeptics\endtitle
\author Saharon Shelah\endauthor
\address Mathematics Department, The Hebrew University, 
Jerusalem, Israel\endaddress
\curraddr{\rm Mathematics Department, Hill Center, 
Busch Campus, Rutgers--The
State University, New
Brunswick, New Jersey 08903}\endcurraddr
\date October 1, 1990 and, in revised form, May 27, 
1991\enddate
\thanks Research partially supported by the BSF; Pub. no. 
400A\endthanks
\endtopmatter

\document

When modern set theory is applied to conventional 
mathematical
problems, it has a disconcerting tendency to produce 
independence
results rather than theorems in the usual sense. The 
resulting
preoccupation with ``consistency'' rather than ``truth'' 
may be felt
to give the subject an air of unreality.  Even elementary 
questions
about the basic arithmetical
operations of exponentiation in the context of infinite 
cardinalities,
like the value of $\cvt$, cannot be settled on the basis 
of the usual
axioms of set theory (ZFC).
 
Although much can be said in favor of such independence
results, rather than undertaking to challenge such
prejudices, we have a more modest goal; we wish to  point 
out an area of
contemporary set theory in which theorems  are abundant, 
although the
conventional wisdom  views the
subject as dominated by independence results, namely, 
cardinal
arithmetic.
 
To see the subject in this light it will be necessary to 
carry out a
substantial shift in our point of view.  To make a very 
rough analogy
with another generalization of ordinary arithmetic, the 
natural
response to the loss of unique
factorization caused by moving from $\Bbb{Z}$ to other 
rings of algebraic
integers is to compensate by changing the definitions, 
rescuing the
theorems.  Similarly, after shifting the emphasis in 
cardinal arithmetic
from the usual notion of exponentiation to a somewhat more 
subtle
variant, a substantial body of results is uncovered that 
leads to new
theorems in cardinal arithmetic and has applications in 
other areas as
well.  The first  shift is from cardinal exponentiation to 
the more
general notion of an infinite product of infinite 
cardinals; the
second shift is from cardinality to cofinality; and the 
final shift is
from true cofinality to potential cofinality (pcf).  The 
first
shift is quite minor and will be explained in \S1.
The main shift in viewpoint will be
presented in \S4 after a review of basics in \S1, a brief 
look at
history in \S2, and some personal history in \S3.
The main results on pcf are presented in \S5.
Applications to cardinal arithmetic are described in \S6. 
The  limitations on
independence proofs are discussed in \S7, and in \S8 we 
discuss
the status of two axioms that arise in the new setting.
  Applications to other areas are found in \S9.
 
The following result is a typical application of the theory.
 
\proclaim{Theorem A} If $2^{\aleph_n}<\Ao$ for all $n$ then
$2^{\aleph_\omega}<\aleph_{\omega_4}$.
\endproclaim
 
 The subscript 4 occurring here is admittedly very 
strange.  Our thesis
is that the theorem cannot really be understood in the 
framework of
conventional cardinal arithmetic, but that it makes 
excellent sense as
a theorem on pcf.  Another way of putting the matter is 
that the
theory of cardinal arithmetic involves two quite different 
aspects,
one of which is totally independent of the\ usual axioms of
set theory, while the other is
quite amenable to investigation on the basis of ZFC. Since 
the usual
approach to cardinal arithmetic mixes these two aspects, the
independence results mask the theorems and the occasional 
result
that survives this interference then looks quite surprising.
 
Of course, the most severe skeptics will even deny the 
mathematical
content of Cantor's theorem ($\cvt >\aleph_0$). To these 
we have
nothing to say at all, beyond a reasonable request that 
they refrain
from using the countable additivity of Lebesgue measure.
 
Most of the results reported here were obtained in the 
past three
years and are expected to appear in a projected volume to 
be published
by Oxford University Press.
 
\heading{1.} Basics\endheading
 
The reader is assumed to be comfortable with the axiom 
system ZFC or
an equivalent formulation of axiomatic set theory 
including the Axiom
of Choice, though knowing naive set theory is enough for 
reading
almost everything .  In particular we have the notions of 
ordinal,
cardinal, and cardinality $|A|$ of a set $A$, the 
identification of
ordinals with sets of ordinals, of cardinals with 
``initial'' ordinals,
and hence also with sets of ordinals. The cofinality 
$\cof\alpha$ of an
ordinal $\alpha$ is min$\{ |A| : A$ is an unbounded subset 
of $\alpha\}$;
we call an infinite cardinal $\lambda$ regular if 
$\cof\lambda =
\lambda$, otherwise we call it singular.
 (Note $\cof\aleph_{\alpha+1}=\aleph_{\alpha+1}$, 
$\cof\aleph_0=\aleph_0$,
 $\cf\aleph_\omega=\aleph_0$, and for
 limit ordinal $\delta$, $\cof\aleph_\delta=\cf\delta)$.
We write
$\alpha,\beta,\gamma,\delta$ for ordinals, with $\delta$ 
typically a
limit ordinal, and $\kappa, \lambda,\mu$ for cardinals.
$\aleph_\alpha$ is the $\alpha$\<th infinite cardinal . 
The first
cardinal above $\lambda$ is denoted $\lambda^+$
$(\aleph_\alpha^+=\aleph_{\alpha+1})$.
 
The product of a set of cardinals is the cardinality of 
their
cartesian product (with each cardinal thought of as a set of
ordinals). Exponentiation is treated as a special case of 
the infinite
product.  We recall that for any ordinal $\alpha$, $\cof 
\alpha$ is a
regular cardinal,
and that $\aleph_\alpha$ is singular iff $\alpha$ is
a limit ordinal with $\cof \alpha <\aleph_\alpha$.
A cardinal $\lambda$ is a limit cardinal
[strong limit cardinal] if
$\mu < \lambda$  implies  $\mu ^+ < \lambda
[ 2^\mu < \lambda ]$.
  Throughout our
discussion the single most interesting limit
cardinal will be
$\aleph_\omega$, as was already illustrated in the 
introduction.
 
As has been known from time immemorial, addition and 
multiplication of
two cardinals trivializes when at least one of them is 
infinite, so
the theory of cardinal arithmetic begins with cardinal 
exponentiation,
and more generally with infinite products. In particular, 
the function
$\lambda\mapsto2^\lambda$ satisfies the following two 
classical laws
(the first being entirely trivial):
$$\hbox{If $\mu\le\lambda$ then $2^\mu\le 
2^\lambda$.}\leqno(1)$$
$$\cof 2^\lambda >\lambda.\leqno(2)$$
 
The most basic problem in (conventional) cardinal 
arithmetic would be
whether this function obeys other laws, and for this one 
quickly comes
to consider the behavior of $2^{\aleph_\delta}$ for 
$\delta$ a limit
ordinal (as the historical discussion in \S2 explains).  
In this case
one has various relations of the type
$$2^{\aleph_\delta}=\prod_{\alpha<\delta} 
2^{\aleph_\alpha}.\tag3a$$
In particular, if $\aleph_\delta$ is a strong limit 
cardinal -- that
is, if $2^{\aleph_\alpha}<\aleph_\delta$ for all $\alpha 
<\delta$ --
we have
$$2^{\aleph_\delta}=\aleph_\delta^{\cof \delta}.\tag3b$$
One of the difficulties in the study of cardinal 
arithmetic is the
preoccupation with $2^\lambda$, or more generally with 
exponentiation
of a small base to a large exponent; the reverse situation 
is
considerably more manageable, and a preoccupation with 
strong limit
cardinals is merely an attempt to trade in one problem for 
the other
via a relation of type (3b).
 Maybe this preoccupation is a vestige
of the Generalized Continuum Hypothesis.
 
We return to the first shift in viewpoint discussed in the
introduction.  If $\aleph_\omega$ is a strong limit 
cardinal then
$2^{\Ao}=\aleph_\omega{}^{\aleph_0}$, by (3b).  In the 
sample
Theorem A given in the introduction, we explicitly 
included the
hypothesis that $\Ao$ is a strong limit cardinal. In fact, 
the
theorem is cleaner without it:
 
\proclaim{Theorem A$'$} $\Ao^\Az <\max 
\{\aleph_{\omega_4},(\cvt)^+\}$.
\endproclaim
 
This is clearly a more ``robust'' statement than the 
original
formulation, in that fewer extraneous considerations are 
involved.
Unfortunately the meaning of the statement will still depend
on the value of $\cvt$, the dependence being trivial if 
$2^{\aleph_0}>\aleph
_\omega$, and we will have to work harder, beginning in
\S4,  to eliminate this feature.
 
 
The following ideas belong to modern cardinal arithmetic, 
though
absent from the classical theory.  There is a natural 
topology on the
class of all ordinals, in which a limit ordinal $\alpha$ 
is a limit
point of any set $X$ of ordinals for which 
$X\intersect\alpha$ is
unbounded below $\alpha$.  If $\alpha$ is an ordinal of 
uncountable
cofinality, $\FF_\alpha$ will denote the filter generated 
by the
closed unbounded subsets of $\alpha$.  This filter is 
closed under
countable intersections.
 
A {\it filter\/} on a set $I$ is a collection of subsets 
of $I$ closed
upward (with respect to inclusion) and closed under 
intersection.
We require filters to be {\it nontrivial\/}, that is, the 
empty set may
not be in the filter.
 
If $\{X_i:i\in I\}$ is an indexed family of sets and $\FF$ 
is any
filter on $I$, then the {\it reduced product\/} $\prod_i 
X_i/\FF$ is the
set of functions $f\in \prod_i X_i$ modulo the equivalence 
relation
$=_\FF$ defined by $f=_\FF g$ iff $\{i:f(i)=g(i)\}\in \FF$.
Similarly, if the $X_i$ are linearly ordered, then 
$\prod_i X_i$ is
partially ordered by pointwise comparison and $\prod_i 
X_i/\FF$ is
partially ordered by pointwise comparison modulo $\FF$.
 
A maximal (nontrivial) filter is also called an 
ultrafilter and the
 corresponding
reduced products are called ultraproducts.  For our 
purposes, the main
point is that an ultraproduct of linearly ordered sets is 
again
linearly ordered and not just partially ordered.
 
We will touch lightly on issues connected with large 
cardinals and
inner models.  The reader who is unfamiliar with these 
subjects may
ignore these remarks.  On the other hand the following 
comments may be
sufficient by way of background. There are various axioms 
concerning
the existence of ``large'' cardinals (bearing names like: 
strongly
inaccessible, measurable, supercompact, huge), which are 
easily seen to
be unprovable on the basis of ZFC; even their consistency 
is known to
be unprovable on the basis of ZFC, though these axioms are 
generally
thought to be consistent.  Some consistency results have 
been obtained
only on the basis of the assumed consistency of large 
cardinal axioms,
and the  application of ``inner models theory'' is a 
method to prove that such
consistency results require large cardinals.  As our main 
concern is
with provable theorems, this takes us rather far afield, 
but there is
a constant interaction between the search for new theorems 
and the
limitations imposed by independence results, which one 
cannot ignore
in practice.  In any case, when we refer to large 
cardinals, it is
understood in that sense.  When we wish to refer to 
cardinals
that are large in a less problematic sense (bigger than some
particular cardinal) we will refer to them as ``moderately 
large.''
 
\heading{2.} History\endheading
 
Until 1974, the classical monotonicity and cofinality 
restrictions  (1),
(2) given
 above exhausted the known
properties of the function $(\lambda\mapsto 2^\lambda)$. 
G\"odel proved
the consistency (with ZFC) of the generalized continuum 
hypothesis:
$2^\Aa=\aleph_{\alpha+1}$. In 1963 Cohen introduced the 
notion of
forcing, setting off the wave of independence results that 
continues
to this day, and used it to prove that $2^\Az$ can be any 
cardinal of
uncountable cofinality. Easton adapted this method to show 
that for
any appropriate
function $f(\lambda)$ satisfying the monotonicity and 
cofinality
restrictions, it is consistent to assume 
$2^\lambda=f(\lambda)$ for all
{\it regular\/} $\lambda$. In Easton's constructions, for 
$\lambda$ singular
 $2^\lambda$ will
always be the least value consistent with (1), (2), and 
the values of
$2^\mu$ for $\mu$ regular. For example, if Easton made 
$\Ao$ a
strong limit cardinal, then he made $2^\Ao=\aleph_{\omega+
1}$.
 
Thus the belief grew that cardinal arithmetic lay outside 
the realm of
mathematical investigation, and to complete the picture it 
apparently
remained only to modify Easton's approach to incorporate 
singular
cardinals.  Under large cardinal hypotheses, progress was 
made first for
 large singular cardinals, and then by Magidor in 1973 for
$\aleph_\omega$, proving for example that one could have 
$\Ao$ a
strong limit and $2^\Ao=\aleph_{\omega+2}$.
 
So it came as a great surprise in 1974 when Silver 
produced a new
theorem of cardinal arithmetic: if $2^{\Aa}=\aleph_{\alpha+
1}$ for all
countable $\alpha$, then 
$2^{\aleph_{\omega_1}}=\aleph_{\omega_1+1}$.
At this point we will leave the later independence results 
aside
(Magidor, Foreman, Woodin, Gitik, Cummings, and the 
present author), as well as
the complementary work on inner models and consistency 
strength
(Jensen, Devlin, Dodd, Mitchell, Gitik), and earlier works 
(Scott,
Solovay and Magidor) and concentrate on theorems provable 
in ZFC.  The
next section concentrates on singular cardinals of 
uncountable
cofinality; one can skip over this development
and continue in \S4.
 
\heading{3.} $\lambda^{\aleph_1}$\endheading
 
In what follows, $\aleph_1$ can be replaced by any 
uncountable regular
cardinal, but not (in this section) by $\omega$.
 
Answering a question of Silver,
 Galvin and Hajnal proved
 
\proclaim{Theorem \rm[GH]} Suppose that $\delta$ is a 
limit ordinal of
cofinality $\aleph_1$, and $\Aa^{\aleph_1}<\Ad$ for 
$\alpha<\delta$.
Then $\Ad^{\aleph_1}<\aleph_{(2^{\mid \delta\mid })^+}$.
\endproclaim
 
More precisely, if we define $f:\delta\to\delta$ by
$\Aa^{\aleph_1}=\aleph_{f(\alpha)}$, the theorem states
that $\Ad^{\aleph_1}<\aleph_{\nr f\nr }$, where $\nr  f\nr 
$ is defined
inductively as $$\sup \{\nr g\nr +1: g < f\, 
\roman{mod}\,\FF_\delta\}.$$  (For
notation, cf. \S1.) This definition turns out to make 
sense because
the filter $\FF_\delta$ is closed under countable 
intersections.
Following Silver, this line was developed
 in parallel and subsequent works
 by Baumgartner, Galvin and Hajnal, Jech and Prikry, 
Magidor, and
the author.
 
When I became interested in the subject, I saw a great 
deal of
activity and suspected I had come into the game too late; 
shortly
thereafter I seemed to be the only one still interested in 
getting
theorems in ZFC.
I believed that the following thesis would be fruitful.
 
\proclaim{Thesis} For $\lambda$ of cofinality $\aleph_1$ 
satisfying
$$\mu^{\aleph_1}<\lambda \hbox{ for all } \mu < 
\lambda\leqno(*)$$
 if $\lambda$ is
small in any sense then $\lambda^{\aleph_1}$ is small in a 
related
sense.
\endproclaim
 
In practice this means one writes $\lambda=F_1(\alpha)$ 
for $F_1$ some
natural function, and proves a bound 
$\lambda^{\aleph_1}<F_2(\alpha)$
where $F_2(\alpha)=F_1((|\alpha| ^{\aleph_1})^+)$ or some 
quite similar
function.
 
For example,
 
\proclaim{Theorem \rm[Sh 111]} If 
$\lambda>2^{2^{\aleph_1}}$ satisfies $(*)$ and is
 below
the first regular uncountable limit
 cardinal, then so is $\lambda^{\aleph_1}$.
 We also get results in this vain for
 $F_1(\alpha)=$ the $\alpha$\<th fix point
 \RM(i.e., cardinal $\lambda$ equal to $\aleph_\lambda)$.
\endproclaim
 
This uses a rank function similar to Galvin/Hajnal's, with 
respect to
more (normal) filters, which we show is well defined 
unless Jensen's work
trivializes the problem.
 
We got a similar bound in [Sh 256] when $\lambda$ is the 
first fixed
point of order $\omega$ and cofinality $\aleph_1$ in the 
enumeration
of the cardinals, solving a problem raised by Hajnal 
following [Sh
111]. There the ranks were with respect to more
complicated objects than normal filters, and in [Sh 333] 
similar
bounds are obtained for functions defined inductively.
We also prove that if the problem is not trivial then
if we collapse $2^{2^{\aleph_1}}$ there is an ultrapower 
of the
old universe in which for all regular 
$\lambda>2^{2^{\aleph_1}}$ there
is a $\lambda$-like element in the ultrapower (in [Sh 111] 
this was
done for each $\lambda$ seperately). In [Sh 386]
$\lambda>2^{2^{\aleph_1}}$ was replaced in the theorem 
above by
$\lambda > 2^{\aleph_1}$. A posteriori the line [Sc,
So, Si, GH, Sh 111, Sh 256, Sh 386] is quite straight.
The rest follows a different line.

\heading{4.}  Possible  cofinalities\endheading
 
Although Cohen and Easton showed us that powers of regular 
cardinals
are easily manipulated, we learned from inner model theory 
that this
is not the case for powers of singular cardinals.  
Similarly,
manipulating $\lambda^\kappa$ for $\lambda >2^\kappa$ is 
much harder
than manipulating $2^\kappa$; the same applies to products 
of
relatively few moderately large cardinals.  There were 
indications
[Sh b, Chapter XIII, \S\S5, 6] that cofinalities are at 
work behind
the scenes.  At a certain point we began to feel that we 
could split
off
the independence results from the hard core of truth by
 shifting the focus.
 
Let $\bfa =\cardset$ be an indexed  set of regular 
cardinals with each
$\lambda_i$ greater than $|I|$.  Where cardinal arithmetic 
is
concerned with the cardinality of $\prod \bfa=\prod_i 
\lambda_i$, we
will be  concerned with the following cardinal invariants 
involving more of the
 structure
of the product.
 
\subheading{\num{4.1.} Definition}
1. A cardinal $\lambda$ is a {\it possible cofinality\/} of
$\prod \bfa$ if there is an ultrafilter $\FF$ on $\bfa$ 
for which the
cofinality of $\prod \bfa/\FF$ is $\lambda$. (Recall that 
these
ultraproducts are linearly ordered: \S1.)

2. $\pcf \bfa$ is the set of all
possible cofinalities of
$\prod \bfa$.
 
We first used pcf in [Sh 68] in a more
structural context, to
construct Jo\'nsson algebras (algebras of a
given cardinality with no
proper subalgebra of the same cardinality).
In [Sh b, Chapter XIII,
\S\S5,6] we obtained results under the
more restrictive hypothesis $\lambda_i>2^{\mid \bfa\mid }$ 
bearing on cardinal
arithmetic. An instance of the main theorem there is the 
following
(note that [GH] did not give information on cardinals with 
countable
cofinalities, Theorem A is a significant improvement 
though formally
they are incomparable).
 
\proclaim{\num{4.2.} Theorem B$'$} If $2^\Az<\Ao$ then 
$\Ao^\Az<\aleph_{(\cvt)^+}$.
\endproclaim
 
Since there are $2^{2^{\mid I\mid }}$ ultrafilters on $I$, 
$\pcf \bfa$ could
be quite large a priori; this would restrict heavily its
applications to cardinal arithmetic.  Fortunately there 
are various
uniformities present that lead to a useful structure 
theory for pcf.
 
\proclaim{\num{4.3.} Main Theorem} 
Let $\bfa=(\lambda_i)_{i\in I}$ be an indexed set of
 distinct regular cardinals, with
each $\lambda\in \bfa$ greater than $|\bfa|$. Then\RM:
\roster
\item"1." $\pcf \bfa$ contains at most $2^{\mid \bfa\mid 
}$ cardinals\RM;
\item"2." $\pcf \bfa$ has a largest element $\max \pcf 
\bfa$\RM;
\item"3." $\cof\prod\bfa=\max\pcf\bfa${\rm (}see remark 
below\RM)\RM;
\item"4." For each $\lambda\in \pcf \bfa$ there is a subset
$\bfb_\lambda$ of $\bfa$ such that
\itemitem"a." $\lambda = \max\pcf \bfb_\lambda$, and
 
\itemitem"b." $\lambda\notin \pcf (\bfa-\bfb_\lambda)$\RM;
\item"5." If $\JJ_\lambda$ is the ideal on $I$ generated 
by the sets
$\bfb_\mu$ for $\mu<\lambda$, then
for each $\lambda\in \pcf \bfa$ there are functions 
$f^{\lambda}_i$
$(i<\lambda)$ such that
\itemitem"a." for $i<j$ we have $f^\lambda_i <
      f^\lambda_j \mod \JJ^\lambda$\RM;
\itemitem"b." for any $f\in \prod\bfa$ and 
$\lambda\in\pcf\bfa$ there
is some $i<\lambda$ such that $f<f^\lambda_i \mod
\langle\JJ_\lambda,(\bfa-\bfb_\lambda)\rangle$.
\endroster

\endproclaim
 
\rem{Remarks} 
The meaning of clause (3) is that for 
$\lambda=\max\pcf\bfa$, there
is a subset $P$ of $\prod\bfa$ of cardinality $\lambda$ 
that is cofinal
in the sense that every function in $\prod \bfa$ is 
dominated
pointwise by a function in $P$.  For example, the 
cofinality of
$\omega\times\omega_1$ is $\omega_1$. This does not  mean 
that
there is a pointwise nondecreasing sequence of length 
$\omega_1$ that
dominates every element of $\omega\times\omega_1$; there 
is no such
sequence. When a product actually contains a cofinal 
pointwise
nondecreasing sequence, we say its cofinality is {\it 
true\/}, and we
write  $\tcf$ for the cofinality when it is true.
\endrem

By clause (5) $\prod\bfb_\lambda/\JJ_\lambda$ has true 
cofinality $\lambda$.
 
The following structural principle is also of great 
practical
importance.
 Note that we do not know whether
   $\mid \pcf\bfa \mid   \le \mid \bfa|$, 
the following still says that $ \pcf\bfa$ is ``small,''
``local,'' has small density character.
 
\proclaim{\num{4.4.} Localization Theorem} Let  $\bfa$ be 
a set of $\kappa$ distinct
 regular cardinals with $\lambda>\kappa$ for all 
$\lambda\in \bfa$\RM;
and suppose $\bfb\includedin \pcf \bfa$ with $\lambda>\mid 
\bfb\mid $ for all
$\lambda$ in $\bfb$. If $\mu\in \pcf \bfb$ then 
$\mu\in\pcf\bfa$, and
for some $\bfc\includedin\bfb$ of cardinality at most 
$\kappa$ we have
$\lambda\in \pcf\bfc$.
\endproclaim
 
Thus $\pcf$ defines a canonical closure operation on sets 
of regular
cardinals with some good properties.
\heading{5.} Pseudopowers\endheading
 
For $\cof \lambda\le \kappa<\lambda$ we define the {\it 
pseudopower\/}
$\pp_\kappa(\lambda)$ as follows.
 
\subheading{\num{5.1.} Definition}
1. $\pp_\kappa(\lambda)$ is the supremum of the 
cofinalities of
the ultraproducts
$\prod\bfa/\FF$ associated with a set of at most $\kappa$ 
regular
cardinals below $\lambda$ and an ultrafilter $\FF$ on $\bfa$
containing no bounded set bounded below $\lambda$.

2. $\pp(\lambda)$ is $\pp_{\cof\lambda}(\lambda)$.
 
If we have a model of set theory with a very interesting 
(read:
bizarre) cardinal arithmetic, say $2^{\Ao}=\aleph_{\omega+
2}$, and
we adjust $2^{\Az}$ by Cohen's method, putting in 
$\aleph_{\omega+3}$
Cohen reals, there is then no nontrivial operation of 
cardinal arithmetic that
will yield the result $\aleph_{\omega+2}$.  Not that the 
original
phenomena have been erased; they are simply drowned out by 
static. The
operation $\pp_\kappa(\lambda)$ is more robust. It will be 
convenient
to call $\lambda$ $\kappa$-inaccessible if 
$\mu^\kappa<\lambda$ for
$\mu<\lambda$.
 
\proclaim{\num{5.2.} Theorem}
\RM1. $\pp_\kappa(\lambda)\le \lambda^\kappa$\RM;

\RM2. If $\cof \lambda\le \kappa<\lambda<\aleph_\lambda$ and
$\lambda$ is $\kappa$-inaccessible
\RM(i.e., $(\forall\mu<\lambda)[\mu^\kappa<\lambda])$, 
{\it then}
$\pp_\kappa(\lambda)=\lambda^\kappa$.
\endproclaim
 
In this sense, $\pp_\kappa(\lambda)$ is an antistatic 
device.  If $\lambda$
is $\kappa$-inaccessible in one
model of set theory and this condition is then destroyed 
in an extension in a
 reasonably subdued
manner (i.e., when we make $2^\kappa$ large), 
$\pp_\kappa(\lambda)$ continues to
reflect the earlier value
of $\lambda^\kappa$.  On another level, one can prove 
results
about this operation by induction on cardinals, which is 
not possible
with the less robust notions. (The restriction
$\lambda<\aleph_\lambda$ is certainly convenient, but we 
will discuss
its removal below.)
 
Let $\PP_\kappa(\lambda)$ be the set of cofinalities whose 
supremum was taken
to get $\pp_\kappa(\lambda)$.  This turns out to have the 
simplest possible
structure.
 
\proclaim{\num{5.3.} Convexity Theorem} If 
$\kappa\in[\cof\lambda,\lambda)$ then
$\PP_\kappa(\lambda)$ is an interval in the set of regular 
cardinals with
minimum element $\lambda^+$.
\endproclaim
 
We can also give a cleaner description of the way a 
cardinal enters
$\PP_\kappa(\lambda)$ in some cases.
 
\proclaim{\num{5.4.} First Representation Theorem}
Suppose $\Az<\cof\lambda<\kappa<\lambda$ and for all $\mu\in
(\kappa,\lambda)$\RM; if $\cof \mu\le \cof\lambda$ then 
$\pp(\mu)<\lambda$.
Let $\FF_0$ be the filter on $\cof \lambda$ generated by 
the complements
of the bounded sets.
Then for any regular cardinal $\lambda^*\in\PP(\lambda)$,
there is an increasing sequence $(\lambda_i)_{i<\cof 
\lambda}$ with
limit $\lambda$ such
that the product $\prod \lambda_i/\FF_0$ has true cofinality
$\lambda^*$.
\endproclaim
 
In the case $\lambda^*=\lambda^+$ there is a better 
representation.
 
\proclaim{\num{5.5.} Second Representation Theorem}
If $\lambda$ is singular of uncountable cofinality, 
$\FF_0$ is the
filter generated by the complements of the bounded subsets 
of $\cof
\lambda$,  and
$(\lambda_i)_{i<\cof \lambda}$ is increasing, continuous  
and cofinal in $\lambda$,
 then there is a
closed unbounded subset $C$ of $\cof \lambda$ such that 
$\prod_C
\lambda_i^+/\langle\FF_0,\roman{cof}\,\lambda-C\rangle$ 
has true cofinality $\lambda^+$.
\endproclaim 
 
\rem{Discussion}
Why do we suggest $PP_\kappa(\lambda)$
as a replacement to $\lambda^\kappa$?
 Maybe we had better reconsider what we are doing.
$\lambda^\kappa$ is a measure of the size of
the family of subset $\lambda$ of cardinality $\le 
\kappa$, i.e., its
cardinality. Remember\ that when $\lambda\le\kappa$ the 
independence results
do not leave us much to say.
When $\kappa<\lambda$
we shall present various natural measures of this set below.
It is useful to prove
that various such measures are equal, as we can then look 
at them as various
characterizations of the same number each useful in 
suitable circumstances.
We report reasonable success in this direction below.
\endrem
\subheading{\num{5.6.} Definition}

For $\kappa\in[\cof\lambda\!,\lambda)$ let the 
$\kappa$-covering number
for $\lambda$\!, $\cov(\lambda\!, \kappa)\!,$ be the 
minimal cardinality of a family of
 subsets of
$\lambda$, each of size less than $\lambda$, such that 
every subset of
$\lambda$ of size $\kappa$ is contained in one of the 
specified sets.
 Note that, as is well known, for $\kappa<\lambda$,
$\lambda^\kappa=2^\kappa\litspace+\litspace$``the 
$\kappa$-covering of $\lambda$.''
 
\proclaim{\num{5.7.} Theorem} If 
$\kappa\in[\cof\lambda,\lambda)$ and
$\lambda<\aleph_\lambda$ then $\pp_\kappa(\lambda)$ is the
$\kappa$-covering number for $\lambda$.
\endproclaim
 
To remove the restriction that $\lambda<\aleph_\lambda$, 
define the
{\it weak\/} $\kappa$-covering number for $\lambda$ as the 
minimal
cardinality of a family of subsets of $\lambda$,
each of size less than $\lambda$, such that every subset of
$\lambda$ of size $\kappa$ is contained in the union of 
countably many
of the specified sets.
 
\proclaim{\num{5.8.} Theorem}
\RM1. For $\kappa$ uncountable, the weak $\kappa$-covering 
number
of $\lambda$ is the supremum of the true cofinalities of 
the reduced products
$\prod \bfa/\FF$ with $\bfa$ a set of at most $\kappa$ 
regular
cardinals below $\lambda$ and $\FF$ a filter on $\bfa$ 
that contains
the complement of every subset of $\bfa$ bounded below 
$\lambda$ and
is closed under countable intersections.

\RM2. If this cardinal is regular, the indicated supremum 
is attained.
\endproclaim
 
Thus our problems are connected with the case of 
cofinality $\omega$.
We give relevant partial information in \S8.
 
We mention
two more invariants that
we can now prove
coincide.
 
\proclaim{\num{5.9.} Theorem}
For $\aleph_0<\kappa\le \lambda$ with $\kappa$ regular the 
following
two cardinals coincide\RM:

\RM1. The minimal cardinality of a family of subsets of 
$\lambda$,
each of size less than $\kappa$, such that any subset of 
$\lambda$ of
cardinality less than $\kappa$ is contained in one of the 
specified
sets.

\RM2. The minimal cardinality of a stationary subset of 
the family
of subsets of $\lambda$ of size less than $\kappa$.
\RM(Stationarity means that for every algebra with set of 
elements $\lambda$
and countably many operations, there is a subalgebra $B$,
$|B|<\kappa$, $B\cap \kappa$ is a ordinal and $B\in S.)$
\endproclaim
 
\heading{6.} Cardinal arithmetic revisited\endheading
 
If $\kappa=\aleph_\alpha$, it will be convenient to call 
$\alpha$ the
{\it index\/} of $\kappa$; otherwise our results tend to 
live entirely
in the land of subscripts.
In [Sh b, Chapter XIII, \S\S5,6] we showed
$$\hbox{The index of $\Ad^{\cof \delta}$ is less than
 $(|\delta|^{\cof\delta})^+$}.$$
In particular, for $\delta=\omega$, the corresponding 
index is below
$(\cvt)^+$ (construed as an ordinal). If $\aleph_\omega$ 
is itself
below $\cvt$ then this contains no information and (our 
usual theme)
inessential modifications of the universe can always make 
this happen.
On the other hand known independence results show that when
$\aleph_\omega$ is a strong limit, the index of 
$(\Ao)^\Az$ can be any
countable successor ordinal, so in principle one could 
hope to prove
that the true bound is $\omega_1$ (which we doubt).  Our 
strongest
result is
 
\proclaim{\num{6.1.} Theorem} $\pp(\Ao)<\aleph_{\omega_4}$.
\endproclaim
 
\proclaim{\num{6.2.} Corollary} 
$\Ao^\Az<\max(\aleph_{\omega_4},(\cvt)^+)$.
\endproclaim
 
The more general formulation is
 
\proclaim{\num{6.3.} Theorem} If $\delta$ is a limit 
ordinal,
$\delta<\aleph_{\alpha+\delta}$ then
 $\pp(\aleph_{\alpha+\delta})<\aleph_{\alpha+\mid 
\delta\mid ^{+4}}$.
\endproclaim
 
For the proof, one looks carefully at the closure 
operation induced on
$\omega_4$ by the $\pcf$ structure on $\aleph_{\omega_4}$ 
(passing to
indices) under the assumption that $\aleph_{\omega_4}$ is
 $\le $
$\pp(\aleph_\omega)$.  There are conflicts between the 
main theorem,
the localization theorem and the second representation 
theorem on the
one hand and  combinatorics of closed unbounded sets inside
$\omega_4$ on the other hand, which yield a contradiction. 
 Ultimately,
the $\pcf$ structure cannot exist on $\omega_4$ because no 
such
closure operation exists. On $\omega_3$, there are two 
questions: does
such a closure operation exist and can it be given by 
$\pcf$?
 
\subheading{Structure  of the proof}
Assuming (toward a contradiction) that $\aleph_{{\omega_4}+
1}$
 belongs to
$\pcf(\Ao)$, for $X\includedin \omega_4$ define $\cl
X=\{i<\omega_4:\aleph_{i+1}\in\pcf(\aleph_{j+1})_{j\in 
X}\}$. By
the Convexity Theorem $\roman{cl}(\omega)=\omega_4$
and by the Localization Theorem, if $i\in \cl X$
then for some countable $Y\includedin X$ we have $i\in \cl 
Y$.
By the Second Representation Theorem, if $\delta<\omega_4$ 
is a limit ordinal
of uncountable cofinality then $\delta$ is the maximal 
element of
$\cl C$ for some closed unbounded subset of $\delta$ as 
well as any
smaller closed unbounded set.  These three properties of 
the closure
operation (alone) eventually lead to a
contradiction.\qed

This may indicate that it is interesting to investigate, 
e.g.,  the set
$\{\aleph_{\alpha+1}:\aleph_{\alpha+1}
\le \roman{pp}\aleph_\omega\}=\pcf\{\aleph_n:n\}$ with the 
relation
``$\lambda\in \pcf\bfb$'' and even the sequence
$\langle\bfb_\lambda:\lambda\in \pcf\{\aleph_n:n\}\rangle$
from 4.3 (for which there are theorems saying ``we can 
choose nicely'').
\heading{7.} For true believers\endheading
 
Naturally our results also give a great deal of 
information regarding
the  types of forcing that are applicable to certain open 
questions.
Any instance of $\lambda\in\pcf\bfa$ normally results from 
some normal
ultrafilter ``at the time when the large cardinals were 
present.''
Some of this information is in the canonical spot in [Sh-b, 
Chapter XIII, \S\S5,6],
 and more
is in [Sh\ 282], for example, if $\pp(\aleph_{\omega_1})$ 
is greater
than $\aleph_{\omega_2}$ then there are many ordinals
$\delta<\omega_2$ of cofinality $\omega$ for which 
$\pp(\Ad)$ is above
$\aleph_{\delta+\omega_1}$.
 
These considerations also shed some light on the problem of
resurrection of supercompactness. By our 9.6(1)
 
\proclaim{\num{7.1.} Corollary} 
If in $V$ we have $\cof(\lambda)<\kappa<\lambda$
and $\pp(\lambda)>\lambda^+$, then there is no universe 
extending $V$ in
which $\lambda^+$ remains a cardinal while $\kappa$ becomes
supercompact or even compact.
\endproclaim
 
Gitik has proved, for example, that if $\cof 
\lambda=\omega$,
$\mu^\Az<\lambda$ for all $\mu<\lambda$, and 
$\pp(\lambda)>\lambda^+$,
then in Mitchell's inner model $o(\lambda)=\lambda^{++}$. 
By his
previous independence results this settles the consistency 
strength of
$\cvt<\aleph_{\omega+1}<\Ao^\Az$ (using [Sh b]), though it 
did not
settle the consistency strength of the full singular 
cardinal
hypothesis, as there are fixed points.  For this purpose 
we have
proved
 
\proclaim{\num{7.2.} Lemma} If $\lambda$ is singular of 
cofinality $\Az$, and if
$\pp(\mu)<\lambda$ for all singular cardinals 
$\mu<\lambda$ of
cofinality $\Az$
 and $pp(\mu ) = \mu^+$ for every large enough
$\mu < \lambda $ of cofinality $\aleph_1$,
 then there is a family of countable subsets of
$\lambda$ containing $\pp(\lambda)$ sets, so that every 
countable
subset of $\lambda$ is contained in one of the specified 
sets.  Thus
if $\mu^\Az<\lambda$ for all $\mu<\lambda$, we find
$\pp(\lambda)=\lambda^\Az$.
\endproclaim
 
Together  this yields an equiconsistency result.
 However  Gitik points out that
this result is intrinsically global. For a localized 
result, we have
 
\proclaim{\num{7.3.} Lemma} If $\lambda$ is singular of 
cofinality $\Az$, and if
$\mu^\Az<\lambda$ for $\mu<\lambda$ and 
$\pp(\lambda)>\lambda^\Az$,
then $\pp(\lambda)$ is quite large\RM; there are at least 
$\aleph_1$
fixed points in the interval from $\lambda^+$ to 
$\pp(\lambda)$.
\endproclaim
 
This still does not settle the behavior of singular 
cardinals of
cofinality $\Az$. The following suggests there are few 
exceptional
values.
 
\proclaim{\num{7.4.} Lemma}
If $\kappa$ is a regular uncountable cardinal and 
$(\lambda_i)_{i\le
\kappa}$ is an increasing continuous sequence of cardinals 
with
$\lambda_i^\kappa<\lambda_\kappa$ for $i<\kappa$, then for 
some closed
unbounded subset $C$ of $\kappa$ we have
$\pp(\lambda_i)=\lambda_i^\kappa$ on $C\cup\{\kappa\}$.
\endproclaim
 
This has application to the construction of ``black boxes.''
\heading{8.} Some hypotheses\endheading
 
\proclaim{\num{8.1.} The Strong Hypothesis}
For all singular cardinals $\lambda$, 
$\pp(\lambda)=\lambda^+$.
\endproclaim
 
This is a replacement for GCH as well as ``$0^{\#}$ does 
not exist''
in some cases. It is weaker than either and is hard to 
change by
forcing, but is consistent with large cardinals (and 
indeed holds
above any compact cardinal [So] while having a more 
combinatorial
character than ``$\neg0^{\#}$''. Of course this will not 
give, e.g., a square
sequence on the successor of a singular cardinal.)
 
\proclaim{\num{8.2.} Lemma} The strong hypothesis implies 
that for any singular
cardinal $\lambda$ and any $\kappa<\lambda$,
 $\lambda^\kappa\le \lambda^++2^\kappa$,
there are $\lambda^+$
subsets of $\lambda$, each of cardinality $\kappa$, such 
that every
subset of $\lambda$ of cardinality $\kappa$ is contained 
in one of the
specified sets\RM; and a Jo\'nsson algebra exists in every 
successor
cardinal.
\endproclaim
 
\proclaim{\num{8.3.} The Weak Hypothesis} 
For any singular cardinal $\lambda$,
there are at most countably many singular cardinals 
$\mu<\lambda$ with
$\pp(\mu)\ge \lambda$.
\endproclaim
 
In my opinion, it is a major problem to determine
whether this follows
from ZFC and is the real problem behind the
determination of the true
bound on $\pp(\Ao)$.
 
\proclaim{\num{8.4.} Lemma}
The weak hypothesis implies\RM:
\roster
\item"$\bullet$"
$\pp(\Ao)<\aleph_{\omega_1}$, and
more generally if $\delta<\aleph_\delta$ then 
$\pp(\aleph_\delta)<\aleph_{\mid \delta\mid ^+}$.
\item"$\bullet$"
$\pcf\bfa$ has cardinality at most $\mid \bfa\mid $.
\item"$\bullet$" $\pp(\lambda)$ has cofinality at least 
$\lambda^+$
for $\lambda$ singular.
\endroster
\endproclaim
 
If we strengthen the weak hypothesis by replacing 
``countable'' by
finite, but only for $\mu$ of cofinality $\Az$,  then
Gitik has proved that the stronger version does not follow 
from ZFC.
\heading{9.} Other applications\endheading
 
If the list of applications below does not contain 
``familiar faces'' the
reader will not lose by
skipping to the concluding remarks.
 
We turn now to results involving more structural 
information, bearing
on almost free abelian groups, partition problems, failure 
of
preservation of chain conditions in Boolean algebras under 
products,
existence of Jo\'nsson algebras, existence of entangled 
linear orders,
equivalently 
narrow Boolean algebras,
and the existence of $L_{\infty,\lambda}$-equivalent 
nonisomorphic models.
 
\proclaim{\num{9.1.} Model theory}
If $\lambda$ has cofinality greater than $\aleph_1$ then 
there are two
$L_{\infty,\lambda}$-equivalent nonisomorphic models of 
cardinality
$\lambda$.
\endproclaim
 
This was known for $\lambda$ regular or 
$\lambda=\lambda^\Az$, and
for strong limit $\lambda$ of uncountable cofinality, but
fails for cofinality $\Az$. There is still a small gap.
 
\proclaim{\num{9.2.} Jo\'nsson Algebras}
There is a Jo\'nsson algebra of cardinality 
$\aleph_{\omega+1}$.
\endproclaim

A Jo\'nsson algebra is an algebra in a finite or countable 
language
that has no proper subalgebra of the same cardinality.
 
This was known previously under the hypothesis 
$2^{\Az}\le\aleph_{\omega+1}$. We
 can
now show that if $\lambda$ is singular and there is no 
Jo\'nsson
algebra of cardinality $\lambda^+$, then $\lambda$ is 
quite large, for
example, $\lambda$ is a limit of weakly inaccessible 
cardinals that do
not admit Jo\'nsson algebras (and there is a Jo\'nsson 
algebra of
inaccessible cardinality $\mu$ if $\mu$ is not Mahlo, or 
even not
$(\mu\times\omega)$-Mahlo).
 
\proclaim{\num{9.3.} Chain conditions}
For any cardinal $\lambda>\aleph_1$, there is a boolean 
algebra that
satisfies
the $\lambda^+$-chain condition while its square does 
not\RM; hence its
Stone-\v Cech compactification has cellularity $\lambda$, 
but its
square has cellularity greater than $\lambda$
\endproclaim

Previous References: [T1, Sh 282] for some singular cases, 
[Sh 280] for regular
above $2^\Az$, [Sh 327] for regular $\lambda\ge \aleph_2$.
 
This really comes from coloring theorems. We give a sample 
of the
latter.
\proclaim{\num{9.4.} Coloring Theorem}
For $\lambda>\aleph_1$ there is a binary symmetric 
function $c$\RM:
$(\lambda^+)^2\to\cof \lambda$ such that for any sequence
$(w_i)_{i<\lambda}$ of pairwise disjoint finite subsets of 
$\lambda^+$
and any $\gamma<\cof \lambda$, there is a pair 
$i<j<\lambda$ with
$c[w_i\times w_j]=\{\gamma\}$.
\endproclaim
 
\proclaim{\num{9.5.} Narrow Boolean algebras}
{\rm1.}
If $\lambda$ is singular and less than $\cvt$, then
there is a $\lambda^+$-narrow Boolean order algebra of 
cardinality
$\lambda^+$ \RM(see below\RM).
The same holds if $\kappa^{+
3}<\cof\lambda<\lambda<2^\kappa$.
 
{\rm2.} The class $\{\lambda$\RM: there is in $\lambda^+$ 
a $\lambda^+$-narrow
Boolean algebra\RM\} is not bounded
\RM(really $(\lambda, \aleph_{\lambda^{+3}+1}]$ is not 
disjoint to it\RM).
\endproclaim
 
\subheading{\num{9.5A.} Definitions}
1. The boolean order algebra associated to a linear order 
$L$
is the boolean algebra of subsets of $L$ generated by the 
closed-open
intervals.

2. A boolean algebra is $\lambda$-narrow if it has no set of
pairwise incomparable elements
of size $\lambda\ (a,b$ are incomparable if $a\not\le b$ and
$b\not\le a)$.

3. A linear order $L$ is {\it entangled\/} if for every 
$n$, for
all choices
of distinct $x_{m,i}$ for $m\le n, i<|L|$, and for any 
subset $w$ of
$\{1,\ldots,n\}$, there are $i<j$ so that for all $m\le n$ 
we have
$x_{m,i}<x_{m,j}$ iff $m\in w$.
 
For $\lambda$ regular, the boolean order algebra 
associated to the
linear order $L$ of cardinality $\lambda$ is 
$\lambda$-narrow if and
only if $I$ is entangled. For background see [ARSh 153, 
BoSh 210, T, Sh
345].
 
The bound $\lambda^{+3}$ in the theorem above is connected 
to \S6. Note
that in addition we can use a subdivision into various 
cases, bearing
in mind that by [BSh 212, T] there is such an algebra of 
cardinality
$\cof 2^\lambda$ if there is a linear order of cardinality 
$2^\lambda$
and density $\lambda$, and usually the absence of such an 
order
enables us to use our current method.
 
\proclaim{\num{9.6.} Almost disjoint sets, almost free 
Abelian groups}
\RM1. If $\lambda$ is singular and $\pp(\lambda)>\lambda^+
$, then there is a
family of $\lambda^+$ subsets of $\lambda$, each of 
cardinality $\cof
\lambda$, such that any $\lambda$ of them admit an 
injective choice
function \RM(called transversal\/\RM). Consequently if 
$\cof \lambda=\Az$ then there is a
$\lambda^+$-free but not free abelian group.

\RM2. The following are equivalent for $\lambda>\mu\ge 
\cvt$.

{\rm a.} There are $\lambda$ subsets of $\mu$ of size 
$\aleph_1$
such that any two have finite intersection.

{\rm b.} For some $n$, there are regular cardinals 
$\lambda_{i,m}$
for  $i<\omega_1$ and $m \le n$, such that for every 
infinite subset
$X$ of $\omega_1$, $\mu\le\max\pcf(\lambda_{i,m}:i\in X, 
m\le n)$.

\RM3. If $\lambda$ is regular,
$2^{<\lambda} < 2^\lambda$, and there is no linear order
of cardinality $2^\lambda$ and density $\lambda$,
then in every regular $\mu \in
[2 ^{<\lambda} , 2^\lambda]$
there is an entangled order.
\endproclaim
 
There are other results on almost disjoint sets and 
$\lambda$-free
abelian groups and a topological question of Gretlis, 
Hajnal, and
Szentmiclossy.
\rem{Concluding Remarks}
The following is  not surprising in view of
 Theorem 6.3, and part of the argument is similar
  but requires, at least from the author, considerably 
more work.
\endrem 
\proclaim{Theorem} If $\delta<\aleph_{\omega_4}$ has 
cofinality
$\aleph_0$, then $\pp(\aleph_\delta)<\aleph_{\omega_4}$, 
and hence the
cofinality of the partially ordered set $\langle S_{\le
\Az}(\aleph_{\omega_4});\includedin\rangle$ is
$\aleph_{\omega_4}$ \RM(where $S_\Az(\lambda)=\{a:a\subseteq
\lambda,\,|a| =\aleph_0\})$. 
\endproclaim
 
%
 
Also the state of Tarski conjecture can be clarified, by 
Jech and
Shelah [JeSh 385]  and,
e.g., if $\roman{pp}(\aleph_{\omega_1}) > 
\aleph_{\omega_2}$ then there is a
Kurepa tree on $\omega_1$.
 
We conclude with some words of the author.
Reflecting on the above it seems that whereas once we knew 
 a considerable
amount about
the uncountable cofinality case and nothing about the 
countable
cofinality case (I mean theorems and not consistency results
or consistency strength results) , now the situation has 
reversed. This
is not accurate---the results of Galvin and Hajnal and those
discussed in \S3 above, are not superceded by the later
results. There does not seem to be a generalization
of ranks to countable cofinality , nor have we suggested
so far anything on fixed points or limits of inaccessibles 
below the contiuum.
Recently we have succeed in getting some results on this 
topic,
  and they will be
reported elsewhere [Sh 420].
See more in Gitik Shelah [GiSh\ 412]
(mainly if $\lambda$ is real valued measurable
$\{2^\sigma:\sigma<\lambda\}$ is finite), [Sh 413]
(more on coloring and Jonsson algebras), [Sh 430]
(on $\cf({\cal S}_{<\aleph_1}(\lambda), \subseteq)\le 
\lambda$
for $\lambda$ real valued measurable, and again of 
existence of trees
 with $\kappa$ nodes and any regular $\lambda\in [\kappa, 
2^\kappa]$
  $\kappa$-branches and improvements of 7.3
and on the smallest values needed for canonization 
theorems).
\Refs 
\REF ARSh 153. U. Abraham, M. Rubin, and S. Shelah, {\it 
On the
consistency of some partition theorems for continuous 
colorings and
the structure of $\aleph_1$-dense real order types\/}, 
Ann. Pure Appl. Logic {\bf 29}
(1985), 123--206.
 
\REF B. J. Baumgartner, {\it Almost disjoint sets, the 
dense set
problem, and the partition calculus\/}, Ann. Pure Appl.
Logic {\bf 9} (1975),
401--439.
 
\REF BP. J. Baumgartner and K. Prikry, {\it On a theorem 
of Silver\/},
Discrete Math. {\bf 14} (1976), 17--22.
 
\REF BoSh 210. R. Bonnet and S. Shelah, {\it Narrow 
boolean algebras\/},
Ann. Pure Appl. Logic {\bf 28} (1985), 1--12.
 
\REF C. P. Cohen, {\it Set theory and the continuum 
hypothesis},
Benjamin, NY, 1966.
 
 \REF Cu. J. Cummings,  {\it   Consistency results in
 cardinal arithmetic\/}, Ph. D. thesis, Cal-Tech,
Pasadena, CA.
 
\REF CW. J. Cummings and H. Woodin,   preprint.
 
\REF DeJ. K. Devlin and R. B. Jensen, {\it Marginalia to a 
theorem of
Silver\/}, Proceedings of the Logic Colloquium, Kiel
1974 (G. H. M\"uller, A. Oberschelp, and K. Potthoff, 
eds.),  Lecture
Notes in Math., vol. 
499, Springer-Verlag, Berlin, 1975, pp. 115--142.
 
\REF DoJ. A. Dodd and R. B. Jensen, {\it The core 
model\/}, Annals Math.
Logic (Ann. Pure Appl. Logic) {\bf 20} (1981), 43--75.
 
\REF FW. M. Foreman and H. Woodin, {\it G.C.H can fail
everywhere\/}, Ann. of Math. (2) {\bf 133} (1991), 1--35.
 
\REF GH. F. Galvin and A. Hajnal, {\it Inequalities for 
cardinal
powers\/}, Ann. of  Math. (2) {\bf 101} (1975), 491--498.
 
\REF GHS. J. Gertlis, A. Hajnal, and Z. Szentmiklossy, 
{\it On the
cardinality of certain Hausdorff spaces\/},
Frolic Memorial Volume, accepted.
 
\REF Gi. M. Gitik, {\it The negation of SCH from 
$o(\kappa)^{++}$},
Ann. Pure Appl. Logic {\bf 43} (1989), 209--234.
 
\REF Gi1. {\Bysame} {\it The strength of the failure of the
 singular cardinal hypothesis\/}, preprint.
 
\REF GM.  M. Gitik and M. Magidor, {\it The singular 
cardinal
 hypothesis revisited\/} (in preparation).
 
\REF GiSh 344. M. Gitik and S. Shelah, {\it On certain 
indestructibility
 of stray cardinals and a question of Hajnal\/}, Arch. 
Math Logic
  {\bf 28} (1989), 35--42.
 
\REF GiSh 412. M. Gitik and S. Shelah,
{\it More on ideals with simple forcing notions\/},
Ann. Pure Appl. Logic (to appear).
 
\REF JP. T. Jech and K. Prikry, {\it Ideals over 
uncountable sets\,\RM:
applications of almost disjoint sets and generic 
ultrapowers\/}, 
Mem. Amer.
Math. Soc.  {\bf 214} (1979).
 
\REF JeSh 385.  T. Jech and S. Shelah, {\it On 
a conjecture of Tarski on products of cardinals\/},
Proc. Amer. Math. Soc. {\bf112} (1991), 1117--1124.
 
\REF Mg. M. Magidor, {\it On the singular cardinals problem
\RM I\/}, Israel J. Math. {\bf 28} (1977), 1--31.
 
\REF Mg1. {\Bysame}  {\it On the singular cardinals problem
\rm II\/}, Ann. of Math. (2) {\bf 106} (1977), 517--547.
 
\REF Mg2. {\Bysame} {\it Chang's conjecture and powers of 
singular
cardinals\/}, J. Symbolic Logic {\bf 42} (1977), 272--276.
 
\REF M. W. Mitchell, {\it The core model for sequences of 
measures \rm I\/},
Math. Proc. Cambridge Philos. Soc. {\bf 95} (1984), 
229--260.
 
\REF Pr. K. Prikry, {\it Changing measurable to accessible 
cardinals\/},
Rozprawy Mat. {\bf LXVIII} (1970), 1--52.
 
\REF Sh 68. S. Shelah, {\it Jo\'nsson algebras in successor
cardinals\/}, Israel J. Math. {\bf 39} (1978), 475--480.
 
\REF Sh 71. {\Bysame} {\it A note on cardinal 
exponentiation\/}, J. Symbolic
Logic {\bf 45} (1980), 56--66.
 
\REF Sh 111. {\Bysame} {\it On powers of singular 
cardinals\/}, Notre Dame
J. Formal Logic {\bf 27} (1986), 263--299.
 
\REF Sh b. {\Bysame} {\it Proper Forcing\/}, Lecture
Notes in Math., vol. 940, Springer-Verlag,
Berlin, 1982.
 
\REF Sh 137. {\Bysame} {\it The singular cardinal problem: 
independence
results\/}, Proceedings of a Symposium on Set Theory,
Cambridge
1978 (A. Mathias, ed.), London Math Soc. Lecture Notes 
Ser., vol. 87,
Cambridge Univ. Press, 
Cambridge and New York, 1983, pp. 116--134.
 
\REF Sh 161. {\Bysame} {\it Incompleteness in regular 
cardinals\/},
Notre Dame J. Formal Logic {\bf 26} (1985), 195--228.
 
\REF Sh 256. {\Bysame} {\it More on powers of singular 
cardinals\/}, Israel
J. Math. {\bf 59} (1987), 263--299.
 
\REF Sh 280. {\Bysame} {\it Strong negative partition 
relations above the
continuum\/}, J. Symbolic Logic {\bf 55} (1990), 21--31.
 
\REF Sh 282. {\Bysame} {\it Successors of singulars, 
productivity of
chain conditions, and cofinalities of reduced products of 
cardinals\/},
Israel J. Math. {\bf 60} (1987), 146--166.
 
\REF Sh 327. {\Bysame} {\it Strong negative partition 
relations below the
continuum\/}, Acta Math. Hungar. in press.
 
\REF Sh 410. {\Bysame} {\it More on cardinal arithmetic\/},
Arch. Math. Logic (to appear).
 
\REF Sh 413. {\Bysame} {\it More Jonsson algebras and 
coloring\/}, preprint.
 
\REF Sh 420. {\Bysame} {\it Advances in cardinal 
arithmetic\/}, to appear.
 
\REF Sh 430. {\Bysame} {\it Further cardinal 
arithmetic\/}, preprint.
 
\REF Sh-g. {\Bysame} {\it Cardinal arithmetic\/}, OUP, (to 
appear).
 
\REF ShW 159. S. Shelah and H. Woodin, {\it Forcing the 
failure of the
CH\/}, J. Symbolic Logic {\bf 49} (1984), 1185--1189.
 
\REF Sc. D. Scott, {\it Measurable cardinals and 
constructible sets\/},
Bull. Acad. Pol. Sci. Ser. Math. Astron. Phys. {\bf 9} 
(1961),
521--524.
 
\REF Si. J. Silver, {\it On the singular cardinal 
problem\/},
Proceedings ICM, Vancouver 1974, vol. I, pp. 265--268,
 
\REF So. R. Solovay, {\it Strongly compact cardinals and 
the GCH\/},
Proceedings of the Tarski Symposium, Berkeley 1971,
Proc. Sympos. Pure Math., vol. xxi, Amer. Math. Soc., 
Providence,
RI, 1974,
pp. 365--372.
 
\REF T. S. Todor\v cevi\`c, {\it Remarks on chain 
conditions in
products\/}, Compositio Math. {\bf 56} (1985), 295--302.
 
\REF T1. {\Bysame} {\it Remarks on cellularity in 
products\/}, Compositio
Math. {\bf 57} (1986), 357--372.
 
\REF W. H. Woodin, {\it The Collected Unwritten Works of 
H. Woodin\/}.

\endRefs
\enddocument